\documentclass{amsart}
\usepackage{graphicx}
\theoremstyle{definition}

\theoremstyle{remark}

\numberwithin{equation}{section}

\begin{document}

\title{On Fubini type theorems for the Riemann integral}

\author{Ilgar Jabbarov, Jeyhun Abdullayev}

\address{Ganja State University}%
\email{ilgar\_js@rambler.ru}%



\maketitle




\textbf{Abstract}

One of the essential questions of the theory of multidimensional integrals concerns the evaluation of integrals taken in given domains. In the simplest case, when integrating over parallelepiped, evaluation can easily be performed by repeated integration. In the case of the Lebesgue integral, the question is easily solvable by Fubini's theorem. In the case of the Riemann integral, the situation is complicated by the difference between Jordan and Lebesgue measures. In this paper, we show that in certain important applications of Riemann integrals, one can establish a modification of the theorem on repeated integration in which Fubini's theorem is as powerful as in the case of the Lebesgue integral.

\maketitle

\section{Introduction}

In problems concerning the evaluation of multidimensional Riemann definite integrals, one needs theorems that reduce them to repeated integration and allow changing the order of integration, if needed. In the simplest case, for double integrals over rectangular domains, the evaluation can be reduced to repeated integration. For example, in [1, p. 261], it is shown that for a continuous function $f$ given in the rectangle $D=[a,b]\times [c,d]$ we have
\[\iint \nolimits _{D}f(x,y)dxdy=\int _{a}^{b}\left(\int _{c}^{d}f(x,y)dy \right)  dx=\int _{c}^{d}\left(\int _{a}^{b}f(x,y)dx \right) dy.\] 

Proof of this statement is based on the definition of definite integral, and passing to the limit under the integral is possible due to the uniform continuity of the function. The same arguments can be applied to multiple integrals over more general cases of domains. 

Let we are given a closed, simple Jordan domain $\Omega $ and another open Jordan domain $\Omega_{0} $ that contains the previous one. Take some continuous function $f(x,y)$ in this wide domain. Then a double integral
\[\int _{\Omega }f(x,y)dxdy \] 
is defined. 

Applying formula (6) ([6, p. 50]) or taking in [2] $r=1$ and $f_{1} (x)=x$, we can rewrite this double integral as a repeated integral as below:
\[\int _{\Omega }f(x,y)dxdy =\int _{a}^{b}dx\int _{\Omega _{x} }f(x,y)dxdy  ,\] 
here by $\Omega _{x} $ is denoted the intersection of the domain $\Omega $ with lines parallel to ordinate axes. Since the domain $\Omega$ is supposed to be simple, this intersection will be some segment $\left[\phi _{1} (x),\phi _{2} (x)\right]$. Therefore, the equality \textit{}
\[\int _{\Omega }f(x,y)dxdy =\int _{a}^{b}dx\int _{\phi _{1} (x)}^{\phi _{2} (x)}f(x,y)dy  \] 
is satisfied.

\section{Auxiliary notes and known results}  
	In studying certain questions, it becomes necessary to use tools such as the change of order of integration under weaker conditions than continuity.  For this reason, the question of whether such an order of integration can be changed, or whether iterated and multiple integrals are coincident, is relevant.  In modern theories, these questions are objects of theorems called Fubini-type theorems (or briefly Fubini's theorems). There were established, in ordinary or improper meaning, various theorems of Fubini-type for the Riemann integrals. In the literature, Fubini-type theorems have been proved under various conditions.

\textbf{Theorem 1}. Let $f(x,y)$ be a function defined in the rectangle 
\[\overline{\boldsymbol{\mathrm{R}}}=\left\{(x,y)|a\le x\le b,\, c\le y\le d\right\}\]
and integrable in this rectangle. Suppose that for all points $x\in [a,\, b]$, with exception of finite number of them, the integral   
\[\int _{c}^{d}f(x,y)dy \] 
exists. Then the iterated integral  
\[\int _{a}^{b}dx \int _{c}^{d}f(x,y)dy \] 
is defined, and the equality below is satisfied:
\[\int _{\Omega }f(x,y)dxdy =\int _{a}^{b}dx \int _{c}^{d}f(x,y)dy .\] 

Changing the order of integration, it is possible to establish the result 
\[\int _{\Omega }f(x,y)dxdy =\int _{c}^{d}dy \int _{a}^{b}f(x,y)dx .\] 
The condition on the existence of a double integral is substantive for the satisfaction of the equality of Theorem 1, [4]. 

In the literature (see, for example, [10]), historically, iterated integrals and their coincidence were considered independently, when we consider integrals depending on parameters, as below:
\[g(y)=\int _{c}^{d}f(x,y)dx. \] 
The question of the integrability of the function $g(y)$ and its calculation reduces to iterated integrals. Their coincidence is an important tool in evaluations. Here, the existence and the coincidence of repeated integrals are questions independent of one another. There are well-known examples that show the independence of one question from another [3].

But, in some conditions, these integrals are closely connected. 

\textbf{Theorem 2}. Let the function $f(x,y)$ be integrable in the rectangle
\[\overline{\boldsymbol{\mathrm{R}}}=\left\{(x,y)|a\le x\le b,\, c\le y\le d\right\}.\] 
Then the integrals 
\[\int _{a}^{b}\int _{c}^{d}f(x,\, y) dxdy,\, \,  \int _{a}^{b}dx\int _{c}^{d}f(x,\, y) dxdy,\, \,  \int _{c}^{d}dy\int _{a}^{b}f(x,\, y) dxdy\,  \] 
are existing and are equal; here, the first integral has an ordinary meaning of a double integral, second and third integrals are taken in following meaning: for every $x$, $a\le x\le b$ we denote $F(x)=\int _{c}^{d}f(x,\, y) dy$ and take the integral of $F(x)$ in the segment  $a\le x\le b$, if it is existing; if for some $x$ the integral does not exist, then $F(x)$ we must define as an arbitrary number placed between lower and upper integrals  $\bar{F}(x)$ and $\underline{F}(x)$.

There are examples showing that if the number $F(x)$ does not belong to this interval, then the statement of the theorem may not be true [8, p.161, problem 1,c]. Moreover, it is established that the function $F(x)$   defined by such way is integrable in Riemann meaning, and the set of points \textit{x,} $a\le x\le b$ for which $\bar{J}(x)\ne \underline{J}(x)$ has zero Lebesgue measure [6, 8].  

Although we formulated Fubini-type theorems for rectangular domains, Theorem 1 can not be applied to Jordan domains of arbitrary view. The main difficulty is connected with the fact that the nature of the intersection of the Jordan domain's boundary with lines parallel to coordinate axes may not have zero measure or may not be measurable in Jordan's sense, which is substantive in some questions. 

In the case of Lebesgue integrals, the theorems of Fubini and Tonelli solve the problem in more general conditions. But in the case of the Riemann integral, the question is complicated and sensitive. In this paper, we show that in some wide type of multiple integrals, it is possible to give such theorems of Fubini-type, which are as powerful as in the Lebesgue case.  Here, one establishes the possibility of such a generalization of  Theorem 1 remaining true even in the case of Jordan domains. As distinct from Theorem 2, partial integrals in this case must be taken in some improper meaning (we shall use the terminology widely used in [9] for one-dimensional integrals). 

Consider a closed bounded Jordan domain $\Omega \subset \boldsymbol{\mathrm{R}}^{2} $. From the definition of a Jordan domain, it follows that the boundary of this domain is possible to cover by the union of quadrates, the sides of which are parallel to coordinate axes, and the sum of squares does not exceed an arbitrary given positive small number $\varepsilon >0$. Denote the union of these quadrates as $E_{\varepsilon } $. Suppose that in the domain  $\Omega $, we are given some integrable function $f(x,y)$. Then we can write the integral over the domain $\Omega $ as an improper integral, and we have the relation:
\[\int _{\Omega }f(x,y)dxdy ={\mathop{\lim }\limits_{\varepsilon \to 0}} \int _{\Omega \backslash E_{\varepsilon } }f(x,y)dxdy .\] 

\section{Main results of the paper}

Our main results are as follows.

\textbf{Theorem 3}. Let $\Omega $ be a bounded closed Jordan domain contained in another open domain $\Omega _{0} $, and the function $f(x,y)$ is continuous in $\Omega _{0} $, being bounded there. Denote
\[a=\inf \{ x|(x,y)\in {\mathop{\Omega }\limits^{\circ }} \} ,\, b=\sup \{ x|(x,y)\in {\mathop{\Omega }\limits^{\circ }} \} ,\] 
where ${\mathop{\Omega }\limits^{\circ }} =\Omega \backslash \partial \Omega $ denotes the inner part of the domain $\Omega $. Then the double integral over the domain $\Omega $ is possible to write as a repeated integral below:
\begin{equation}\label{GrindEQ__1_}
\int _{\Omega }f(x,y)dxdy =\int _{a}^{b}\left(\int _{\Omega _{x} }f(x,y)dy \right) dx;
\end{equation} 
here, in the right hand side, the inner integral, with the variable of integration \textit{y},\textit{ }is taken over the intersection $\Omega _{x} =\{ y|(x,y)\in {\mathop{\Omega }\limits^{\circ }} \} $ in following ``improper meaning'':
\[\int _{\Omega _{x} }f(x,y)dy ={\mathop{\lim }\limits_{\varepsilon \to 0}} \int _{\Omega _{x} \backslash E_{\varepsilon } }f(x,y)dy, \] 
moreover, when $\varepsilon<\varepsilon'$ we must have $E_{\varepsilon}\subset E_{\varepsilon'}$.

Before proving the theorem, let us explain the existence of the integral in ``improper meaning''. Note that, as a covering, we take a finite family of quadrates, each of which contains some point of the set $M$ we cover. We assume that no quadrate of the family contains the others. It is clear that for every quadrate $D_i$, it is possible to cover the intersection $M\cap D_i$ by a new family of quadrates contained in $D_i$ with total measure less than the measure of $D_i$. By this reason, when $\varepsilon<\delta$ we have $E_{\varepsilon}\subset E_{\delta}$. Since the function $f(x, y)$ is bounded, and the domain of integration is increasing as $\varepsilon\to 0$, the limit in the formulation of Theorem 3 exists. 

Note that the expression ``improper meaning'' in the formulation of Theorem 3 does not have the ordinary meaning of improper integral in the Riemann sense, because in the ordinary meaning case, we must take the intersection of the line with the abscissa $x$ with boundary, and the nature of such intersection is not known. This integral exists for all values of $x$.

\textbf{Proof.} Let $\varepsilon >0$ be some positive number. Since the function $f(x,y)$ is uniformly bounded in the open domain $\Omega_0 $, then  there exists a constant $K >0$ for which
\[\left|f(x,y)\right|\le K,\] 
for every pair $(x, y) \in \Omega_0 $. Since the domain $\Omega $ is Jordan, then its boundary can be covered by the union of quadrates, with total measure not exceeding $\varepsilon >0$, contained in the domain $\Omega_0$. Denote this quadrates as $\Delta _{1} ,...,\Delta _{N} $. So we have
\[\partial\Omega \subset \bigcup _{n=1}^{N} \Delta _{n} ;\, \sum _{n=1}^{N}\mu(\Delta _{n} )<\varepsilon  ,\] 
where $\mu$ is a Jordan measure. Denoting the union of quadrates of the covering by $\Delta $, we can write
\[\left|\int _{\Delta }f(x,y)dxdy \right|\le K\mu(\Delta )\le K\varepsilon .\] 
Consider the set ${\mathop{\Omega }\limits^{\circ }} \backslash \Delta $. This difference can be represented as a union of rectangles $D_{1},..., D_{M} $ intersecting only along their boundaries. Since the function $f(x,y)$ is continuous in the domain $\Omega $, then it is integrable in every one of these rectangles, and every integral can be represented as repeated integrals. Let for the fixed $x$, such that the set $\Omega_{x,\varepsilon}=\{ y|(x,y)\in {\mathop{\Omega }\limits^{\circ }} \backslash \Delta \}$ is not empty, denote by $D'_{1},..., D'_{k} $ rectangles which have nonempty intersections with the line parallel to the ordinate axes. We have, in accordance with Theorem 1: 
\[\iint \nolimits _{D'_{i} }f(x,y)dxdy =\int _{a_{i} }^{b_{i} }\left(\int _{c_{i} }^{d_{i} }f(x,y)dy \right) dx;\] 
here $D'_{i} =[a_{i} ,b_{i} ]\times [c_{i} ,d_{i} ]$. Summing up, we obtain:
\[\iint \nolimits _{{\mathop{\Omega }\limits^{\circ }} \backslash \Delta }f(x,y)dxdy=\sum _{i=1}^{k}\iint \nolimits _{D'_{i} }f(x,y)dxdy=   \] 
\[=\sum _{i=1}^{k}\int _{a_{i} }^{b_{i} }\left(\int _{c_{i} }^{d_{i} }f(x,y)dy \right) dx=\int _{a}^{b}\left(\sum _{i,(\exists y)((x,y)\in D'_{i} )}\int _{c_{i} }^{d_{i} }f(x,y)dy  \right)  dx.\] 
Note that when the variable $x$ varies in the segment $[a, b]$ in the inner sum, all of the quadrates of the covering, the union of which coincides with the difference ${\mathop{\Omega }\limits^{\circ }} \backslash \Delta $, take part. So, we can write
\[\left|\iint \nolimits _{\Omega }f(x,y)dxdy -\int _{a}^{b}\left(\int _{\Omega _{x,\varepsilon} }f(x,y)dy \right) dx\right|\le \iint \nolimits _{\Delta }\left|f(x,y)\right| dxdy\le K\varepsilon .\] 
Tending $\varepsilon \to 0$, we obtain:
\[\int _{\Omega }f(x,y)dxdy =\lim_{\varepsilon\to 0}\int _{a}^{b}\left(\int _{\Omega _{x, \varepsilon} }f(x,y)dy \right) dx.\]
Now, we note that, for every $x$, the following limit exists:
\[h(x)=\lim_{\varepsilon\to 0}\int _{\Omega _{x, \varepsilon} }|f(x,y)|dy,\]
by the theorem of Lebesgue on bounded convergence ([6, p. 427]). Therefore, we have proved the relation of Theorem 3 in the Lebesgue sense:
\[\int _{\Omega }f(x,y)dxdy =(L)\int _{a}^{b}h(x) dx.\]

Now, for the completion of the proof, we must prove that on the right-hand side, the integral can be taken in the Riemann sense. Let
\[\varphi(x)=\int _{\Omega _{x} }f(x,y)dy\]
and
\[\varphi_{\varepsilon}(x)=\int _{\Omega _{x, \varepsilon} }f(x,y)dy.\] 
Consider some dissection of the segment $[a, b]$: $a=x_0<x_1<\cdots<x_{n-1}<x_n=b$. Taking integral sums for the functions $\varphi$ and $\varphi_l$, we obtain:
\[\left|\sum _{k=0}^{n-1}\left(\varphi (\xi _{k} )-\varphi _{\varepsilon} (\xi _{k} )\right)\Delta x_{k} \right|\le \sum _{k=0}^{n-1}\left|\left(\varphi (\xi _{k} )-\varphi _{\varepsilon} (\xi _{k} )\right)\right|\Delta x_{k}  \le \] 
\[\le \sum _{k=0}^{n-1}\int _{(\xi_k,y)\in E_{\varepsilon}} \left|f (\xi _{k} ,y)\right|dy \Delta x_{k} \le K\varepsilon;\]  
here $\xi_k\in [x_{k-1}, x_k]$. The sum $\sum _{k=0}^{n-1}\varphi _{\varepsilon} (\xi _{k} )\Delta x_{k} $ tends to the integral 
\[\int_a^b\varphi_{\varepsilon}(x)dx,\]
as $n\to \infty$ ($\max_{k}|x_{k-1}- x_k|\to 0$).

The sum $\sum_{k=0}^{n-1}\phi(\xi_k)\Delta x_k$ is bounded for every $n$. Therefore, there exists a subsequence $n_s$ such that this sum tends to some finite limit $I$. So, passing to the limit in the inequality
\[\left|\sum_{k=0}^{n_s-1}\phi(\xi_k)\Delta x_k-
\sum_{k=0}^{n_s-1}\phi_{\varepsilon}(\xi_k)\Delta x_k\right|\le K\varepsilon\]
first as $s\to\infty$ and after $\varepsilon\to 0$, we obtain
\[I=\lim_{\varepsilon}\lim_{s\to\infty}\sum_{k=0}^{n_s-1}\phi_k(\xi_k)\Delta x_k=
\lim_{\varepsilon\to \infty}\int_a^b\phi_{\varepsilon}(x)dx.\]
Since the limit $I$ does not depend on the subsequence $n_s$, we see that the function $\varphi(x)$ is integrable in the Riemann sense and its value is equal to \eqref{GrindEQ__1_}. Theorem 3 is proved.

Shortage of Theorem 3 consists in supposing the continuity of the function $f(x,y)$. We can omit this condition and prove an analog of Theorem 3 for integrable functions. This is a generalization of the Fubini theorem for the Riemann integral, which has almost the same power as in the Lebesgue integral case. Moreover, partial integrals exist in an improper sense, for all \textit{x}. This is a very important generalization of the Fubini theorem, useful in applications. Using a special case of a partial improper integral in [2], the new notion of improper surface integral is defined. 

For the formulation of the Fubini theorem of the type of Theorem 3, we need some remarks and designations. The following lemma of Lebesgue ([6, p.36]) gives a criterion for integrability in the Riemann sense.

\textbf{Lemma 1. }For the integrability of a bounded, in Jordan measurable closed domain $\Omega $, function $f(x, y)$, it is necessary and sufficient that the set of points of discontinuity of this function has zero Lebesgue measure.

Let the function $f(x,y)$ be integrable in some Jordan domain contained in the rectangle
\[\overline{\boldsymbol{\mathrm{R}}}=\left\{(x,y)|a\le x\le b,\, c\le y\le d\right\}.\] 
We impose the condition that the set of points of discontinuity of this function has zero Jordan measure. It means that for every positive number $\varepsilon >0$, there exist finite family of quadrates $(\Delta _{m} )$ with total measure less than $\varepsilon $:
\[\sum _{m=1}^{M_{\varepsilon} }\mu(\Delta _{m}) <\varepsilon  ,\] 
the union of which contains all of the points of discontinuity; moreover, we can suppose that this family does not contain the quadrate included by any other one. Taking some point $x$ from the interval $[a, \textit{b}]$, consider a section of the rectangle $\bar{\boldsymbol{\mathrm{R}}}$ by a line parallel to the ordinate axes. It is clear, due to the condition imposed above, that every line, being parallel to the ordinate axes, intersects with no more than a finite number of quadrates of covering (if they exist). So, the line, parallel to the ordinate axes, with the abscissa $x$, contains a finite number of intervals, contained by quadrates of the covering. These intervals are disjoined by no more than a finite family of segments of continuity (some of these segments may be reduced to the points). In the segments, the function is continuous and therefore integrable. Since every line parallel to the ordinate axes with abscissa $x$ intersects with sides of a finite number of open quadrates, then all such lines with abscissa from a neighborhood of $x$ will also intersect with sides of the same quadrates. As the point $x$ varies, the union of all closures of these neighborhoods, taken together, covers the segment [a, b]. Let  $a\le x_{0} \le b$ is arbitrary point. Then for some sufficiently small $\delta >0$ every line of a view $x=x_{0} +u,\, \left|u\right|<\delta $ intersects with one and the same quadrates. So, we can find a finite number of closed rectangles that are separated from the rectangle
\[R:[x_{0} -\delta ,x_{0} +\delta ];c\le y\le d\] 
(we take $a$ instead of $x_0-\delta$ when $x_0-\delta\le a$ and $b$ instead of $x_0+\delta$ when $x_0+\delta >b$) by deleting open rectangles along which this strip intersects quadrates. In these rectangles, the function is continuous. So, the integral 
\[\int _{\Omega _{x,\varepsilon} }f(x,y)dy \] 
exists for all $x\in [x_{0} -\delta ,x_{0} +\delta ]$. Let us denote $\Delta =\bigcup _{m}^{}\Delta _{m}  $. By Theorem 3,
\[\int _{R}f(x,y)dxdy =\int _{x_0-\delta}^{x_0+\delta}\left(\int _{R\bigcap \Omega _{x,\varepsilon} }f(x,y)dy \right) dx+\iint \nolimits _{\Delta \bigcap R}f(x,y)dxdy .\] 
Summing over all segments of a view $[x_{0} -\delta ,x_{0} +\delta ]$ which cover the segment $a\le x\le b$, we obtain: 
\begin{equation} \label{GrindEQ__2_} 
\left|\int _{\bar{\boldsymbol{\mathrm{R}}}}f(x,y)dxdy -\int _{a}^{b}\left(\int _{\Omega _{x,\varepsilon} }f(x,y)dy \right) dx\right|\le \iint \nolimits _{\Delta }\left|f(x,y)\right|dxdy\le K\varepsilon  .                      
\end{equation} 

The set $\Omega _{x,\varepsilon} $ is ``increasing'' as $\varepsilon \to 0$. Taking some decreasing sequence of real numbers $\varepsilon _{1} >\varepsilon _{2} >\cdots $ such that $\varepsilon _{n} \to 0$, we have:
\begin{equation} \label{GrindEQ__3_} 
\Omega _{x} =\bigcup _{n}\Omega _{nx}  ;                                                           
\end{equation} 
here $\Omega _{nx} $ means the set $\Omega _{x, \varepsilon} $ defined above for the value $\varepsilon =\varepsilon _{n} $. Passing to the limit as $\varepsilon _{n} \to 0$, we deduce from \eqref{GrindEQ__1_}:
\[\int _{\bar{\boldsymbol{\mathrm{R}}}}f(x,y)dxdy =\int _{a}^{b}\left(\int _{\Omega _{x} }f(x,y)dy \right) dx,\] 
in an improper meaning. 

Now, as above, we put $\varphi(x)=\int _{\Omega _{x} }f(x,y)dy$ and $\varphi_l(x)=\int _{\Omega _{lx} }f(x,y)dy$. Then the relations
\[\left|\sum _{k=0}^{n-1}\left(\varphi (\xi _{k} )-\varphi _{l} (\xi _{k} )\right)\Delta x_{k} \right|\le \sum _{k=0}^{n-1}\left|\left(\varphi (\xi _{k} )-\varphi _{l} (\xi _{k} )\right)\Delta x_{k} \right| \le \] 
\[\le \sum _{k=0}^{n-1}\int \nolimits _{\Delta_l} \left|f (\xi _{k} ,y)\right|dy \Delta x_{k} \le M\varepsilon _{l} \to 0\]  
as $n\to \infty$ first, and $l\to \infty$ after, show that the function $\varphi(x)$ is integrable (here $\Delta_l$ denotes the union of cubes with total measure less than $\varepsilon_l$). So, we have proved the following theorem.

 \textbf{Theorem 4}. Let the function $f(x,y)$ be integrable in some closed Jordan domain $\Omega$  contained in the open rectangle
\[\boldsymbol{\mathrm{R}}=\left\{(x,y)|a< x< b,\, c< y< d\right\}.\] 
Suppose that the set of points of discontinuity has zero Jordan measure. Then the double integral over the domain $\Omega $ is possible to write as a repeated integral below:
\[\int _{\Omega }f(x,y)dxdy =\int _{a}^{b}\left(\int _{\Omega _{x} }f(x,y)dy \right) dx;\] 
here on the right-hand side, the inner integral, with the variable of integration \textit{y},\textit{ }is taken over the set $\Omega _{x} $, defined above by \eqref{GrindEQ__3_}, in ``improper meaning'', and represents an integrable function with respect to $x$. 

 \textbf{Theorem 5. }Let the conditions of Theorem 4 be satisfied in the rectangle
\[\boldsymbol{\mathrm{R}}=\left\{(x,y)|a< x< b,\, c< y< d\right\}.\] 
Suppose that the set of points of discontinuity is closed. Then the double integral over the domain $\Omega $ is possible to write as a repeated integral below:
\[\int _{\Omega }f(x,y)dxdy =\int _{a}^{b}\left(\int _{\Omega _{x} }f(x,y)dy \right) dx;\] 
here on the right-hand side, the inner integral, with respect to \textit{y}, is taken over the set $\Omega _{x,\varepsilon} $, defined above, in an improper sense and represents an integrable function in the Riemann sense.

 \textbf{Proof. }Since the function $f(x,y)$ is integrable, the set of points of discontinuity of this function has zero Lebesgue measure. But this set is closed. For this reason, it has zero Jordan measure. Consequently, for any positive real number $\varepsilon >0$, the set of points of discontinuity may be covered by a finite number of quadrates with total measure less than $\varepsilon $. Then every line parallel to the ordinate axes can intersect these quadrates through no more than a finite number of intervals. Then we can apply the reasoning of the proof for Theorem 4. This remark completes the proof.
 
Theorems 3, 4, and 5 take precedence over Theorem 2 in applications, which is obvious from their formulation. Despite Theorem 2 establishing the Fubini theorem for the Riemann integral in its full generality, for applications, the definition of the function $F(x)$ in that form for the repeated integral is not appropriate. But in formulations of Theorems 3 and 4, this integral has a definite meaning even though it is taken in an improper sense. Despite the example given after the formulation of Theorem 2, formulations of Theorems 3, 4, and 5 reflect the classic Fubini theorem. Consider an example in which, for some values of $x$, the value of the partial integral, that is, the integral 
\[\int _{\Omega _{x} }f(x,y)dy \] 
is not placed in the segment [$\bar{F}(x)$, $\underline{F}(x)$]. For the construction, we shall use examples of functions having a given closed subset as a set of points of discontinuity ([3]). 

\begin{figure}[h]
\centering
\includegraphics[width=0.8\textwidth]{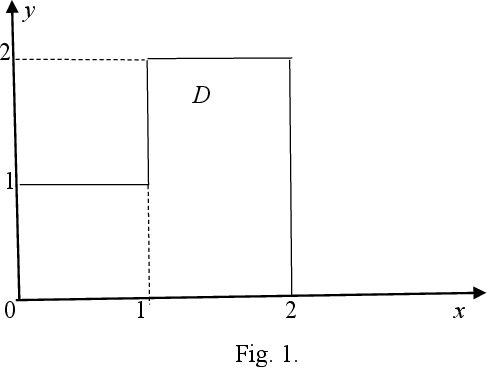}
\end{figure}

Take the domain \textit{D }of a view (Fig.1), which is the union of the rectangle $[0,2]\times [0,1]$ and the quadrate $[1,2]\times [1,2]$. This domain is a Jordan domain. Note that the notion of a Jordan domain is defined in various ways by various authors. For example, in [3], it is defined by simple curves. We accept the definition given in [7, 8]. On the left side of the quadrate, for $x=1$ take some perfect subset of positive measure denoting it as $P$: $P\subset [1,2],\, m(P)>0$. We put 
\[f(x,y)=\left\{\begin{array}{c} {-2,\, \, \, \, \, \, \, \, \, \, \, \, \, \, \, \, \, \, \, \, \, \, \, \, \, \, \, \, \, \, \, \, \, \, \, \, \, \, \, \, \, \, \, \, \, \, \, \, \, \, \, \, \, \, \, \, \, \, \, \, \, for\, x=1,y\in P,} \\ {-1,\, \, \, \, \, \, \, \, \, \, \, \, \, \, \, \, \, \, \, \, \, \, \, \, \, \, \, \, \, \, \, \, \, \, \, \, \, \, \, \, \, \, \, \, \, \, \, \, \, \, \, \, \, \, \, \, \, \, \, \, \, for\, x=1,\, y\notin P,} \\ {0,\, \, \, \, for\, 0\le x\le 2,\, 0\le y\le 1\, and\, 1<x\le 2,\, 1\le y\le 2.} \end{array}\right\} \] 
It is obvious that the function defined by these equalities is integrable and
\[\iint \nolimits _{D}f(x,y)dxdy=0 .\] 
This function is discontinuous at every point (1,\textit{y}), where \textit{y} belongs to \textit{P}. By the theorem of Lebesgue, the partial integral 
\[\int _{0}^{2}f(1,y)dy \] 
does not exist, because the function \textit{f}(1,\textit{y}) has the set of points of discontinuity having positive Lebesgue measure. For lower and upper integrals $\underline{F}(x)$ and $\bar{F}(x)$ we have:
\[\underline{F}(1)=-2\int _{1}^{2}dy =-2;\] 
\[\bar{F}(1)=-\int _{1}^{2}dy =-1.\] 
But at the point \textit{x}=1, one finds
\[F(1) =\int _{0}^{1}f(1,y)dy =0.\] 
Applying Theorem 5, we obtain:
\[0=\iint \nolimits _{D}f(x,y)dxdy=\int _{0}^{2}\left(\int _{\Omega _{x} }f(x,y)dy \right)  dx,\] 
moreover $F(1) \notin [\bar{F}(1), \underline{F}(1)$]. 

For establishing a more general form of Fubini type theorems, we need to introduce some notions. Introduce first the notion of the oscillation of the function at a point. In consent with the definition, the oscillation of the function $f(x, y)$ at the point $(x, y)\in K$ is a quantity:
\[\omega (x, y)=\lim_{\begin {array}{c}{\delta\to 0} \\ {\delta >0} \end{array}}|M_{\delta}-m_{\delta}|;\]
where $M_{\delta}$ and $m_{\delta}$ mean the $\sup$ and the $\inf$ for that function, correspondingly.
In [6, p. 228], the sets $E_{\lambda}, \lambda >0$ were introduced by using of the equality:
\[E_{\lambda}=\left\{(x, y)\in K |\omega (x,y)\ge \lambda \right\}.\]
So, the symbol $E_{\lambda}$ denotes the set of points at which the oscillation of the function is not less than $\lambda$. In [6], the following two statements are proved (theorem 5, p. 288).

\textbf{Lemma 2. } For the continuity of the function $f$ at the point $(x, y)$, it is necessary and sufficient that the oscillation of the function is equal to zero.

\textbf{Lemma 3. } For every $\lambda >0$, the set $E_{\lambda} $ is a closed subset of the set of points of discontinuity.                                                                              

We shall consider below the more general case when the function is integrable, without any conditions on the set of points of discontinuity, as above. When we use Lemma 1, it is possible to cover the set of the points of discontinuity by denumerable family $\Phi$ of open quadrates with total area less than arbitrary given positive number $\varepsilon$. According to Lemma 3, the set $E_{\lambda}$, with $\lambda >0$, has the measure 0, and by compactness, the set $E_{\lambda}$ is possible to cover by some finite subfamily of the family $\Phi$. Denote by $\Sigma_{\varepsilon}$  the union of the all quadrates of the family $\Phi$. Take, as above, an infinite sequence of decreasing positive numbers $\varepsilon_1 >\varepsilon_2 >\cdots$. We can suppose that the sequence of subsets $\Sigma_{\varepsilon_n}$ is decreasing. 

Note that $E_{\lambda}\subset E_{\lambda'}$ if $\lambda' <\lambda$. Fix now arbitrary sequence of positive numbers $\lambda_1 >\lambda_2 >\cdots$. It is clear that the set 
\[\bigcup_{n=1}^{\infty}E_{\lambda_n}\]
coincides with the set of the all points of discontinuity. For every natural $s$, denote by $\Sigma_{\varepsilon, s}$ the finite covering of the set $E_{\lambda_s}$ by quadrates from the family $\Phi$. 

For every fixed value of $x$, we introduce the function:
\[\Gamma_{s,\varepsilon}(x)=\int^{\ast}_{K\setminus\Sigma_{\varepsilon, s}}f(x, y)dy; \]
where on the right hand side, the upper integral stands. When $s$ increases, from the properties of the Darboux sums it follows the inequality $\Gamma_{s,\varepsilon}(x)\geq\Gamma_{s+1,\varepsilon}(x)$; so the limit below exists:
\[\Gamma_{\varepsilon}(x) =\lim_{s\to\infty}\Gamma_{s,\varepsilon}(x).\]
Put now:
\begin{equation} \label{GrindEQ__4_}
\Gamma(x)=\lim_{n\to\infty}\Gamma_{x,\varepsilon_n}.
\end{equation}
This "improper integral" exists for all values of $x$.

The following theorem is a generalization of the theorems established above. 

\textbf{Theorem 6. } Let the function $f(x,y)$ be integrable in he rectangle 
\[\overline{\boldsymbol{\mathrm{R}}}=\left\{(x,y)|a\le x\le b,\, c\le y\le d\right\},\] 
which is included into some open domain.
Then the function $\Gamma(x)$ is integrable in the segment $[a, b]$ and the double integral over $K$ is possible represent as the following iterated integral:
\[\iint _{K}f(x,y)dxdy =\int _{a}^{b}\Gamma(x) dx.\] 

\textbf{Proof.} The double integral under consideration can be represented as follows:
\[\iint_K f(x, y)dxdy=\lim_{\varepsilon\to 0}\iint_{K\setminus \Sigma_{\varepsilon}}f(x, y)dxdy.\]
Really, since the function is bouded, we can denoting by $M$ the maximal value of the modulus of the function $f(x,y)$, write, using Lebesgue integral:
\[|\iint_{\Sigma_{\varepsilon}}f(x, y)dxdy|\le M\varepsilon;\]
so, the limit above is existing.

In [6, p.287], it was proved that for any positive number $\lambda$ there exists such a number $\delta>0$, for which the inequality
\[\sup |f(x, y)-f(x', y')|\le\lambda\]
is satisfied for every pair of tuples $(x, y)$ and $(x', y')$, if they satisfy the inequality $|x-x'|+|y-y'|<\delta$; 
the $\sup$ is taken over all pairs satisfying the above condition. For every natural $s$, we denote the finite covering of the set $E_{\lambda_s}$ by $\Sigma_{\varepsilon, s}$ and put:
\[F_s(x, y)=\left\{\begin{array}{c}{0, if\,\, (x,y)\in\Sigma_{\varepsilon, s}},\\{f(x, y)\,\,\,\, if else.}\end{array}\right\}\]
Every function $F_s(x, y)$ is integrable in $K$, due to Lebesgue criterian. Really, since the set $\Sigma_{\varepsilon, s}$ is closed, all limit points of any sequens from this set belong to this set. Note that the function $F_s(x, y)$ is continuous in the inner part of the set $\Sigma_{\varepsilon, s}$. So, points of discontinuity of this function are the points of discontinuity of the function $f(x, y)$ and the possible points of the boundary of the set $\Sigma_{\varepsilon, s}$. So, the measure of the set of points of discontinuity is equal to zero. We have
\begin{equation} \label{GrindEQ__5_}
\lim_{s\to\infty}\iint_K F_s(x, y)dxdy=\iint_{K\setminus\Sigma_{\varepsilon}}f(x, y)dxdy,
\end{equation}
because the sets $\Sigma_{\varepsilon, s}$ increasingly settles the set $\Sigma_{\varepsilon}$ (note that on the right hand side the Lebesgue integral stands). 

Define in $K$ the dissection $T$, the parameter of which does not exceed the $\delta$ which is defined by $\lambda_s$. Applying the reasoning of [8, p.155], we find:
\[\underline{D}_T(F_s)=\sum_{i=1}^n\sum_{j=1}^m\inf_{\xi_i\eta_j}F_s(\xi_i,\eta_j)
\Delta x_j\Delta y_j\le\]
\[\le\sum_{i=1}^n\inf_{\xi_i}\left(\sum_{j=1}^m\inf_{\eta_j}F_s(\xi_i,\eta_j)\Delta y_j\right)
\Delta x_i\le\]
\[\le\sum_{i=1}^n\inf_{\xi_i}\left(\int_{\ast K\setminus\Sigma_{\varepsilon}} F_s(\xi_j,y)dy\right)\Delta x_i\le\]
\[\le\sum_{i=1}^n\sup_{\xi_i}\left(\int^{\ast}_{K\setminus\Sigma_{\varepsilon}} F_s(\xi_j,y)dy\right)\Delta x_i\le\]
\[\le\sum_{i=1}^n\sup_{\xi_i}\left(\sum_{j=1}^m\sup_{\eta_j}F_s(\xi_i,\eta_j)\Delta y_j\right)
\Delta x_i=\]
\begin{equation} \label{GrindEQ__6_}
=\sum_{i=1}^n\sum_{j=1}^m\sup_{\xi_i\eta_j}F_s(\xi_i,\eta_j)
\Delta x_j\Delta y_j=\overline{D}_T(F_s);
\end{equation}
here the sign $\ast$ below denotes the lower integral, and the sign $\ast$ above denotes the upper integral. It is obvious that
\[|\overline{D}_T(F_s)-\underline{D}_T(F_s)|\le\]
\[\le\sum_{i=1}^n\sum_{j=1}^m
(\sup_{\xi_i\eta_j}F_s(\xi_i, \eta_j)-\inf_{\xi_i\eta_j}F_s(\xi_i, \eta_j))\Delta x_i\Delta y_j\le\lambda_s|K|;\]
here $|K|$ is an area of the rectangle. For every pair $(x, x')$ satisfying the condition $|x-x'|<\delta$, we can write:
\[|\sum_{j=1}^m\inf_{\eta_j}F_s(x, \eta_j)\Delta y_j-\sum_{j=1}^m\sup_{\eta_j}F_s(x', \eta_j)\Delta y_j|\le\]
\begin{equation} \label{GrindEQ__7_}
\le\sum_{j=1}^m|\inf_{\eta_j}F_s(x, \eta_j)-\sup_{\eta_j}F_s(x', \eta_j)|\Delta y_j\le\lambda_s(d-c)\to 0, 
\end{equation}
as $s\to\infty$. 

Applying Lebesgue theorem on monotone convergence for the functions introduced above, we conclude
\[\int_a^b\Gamma_{\varepsilon}(x)dx=\lim_{s\to\infty}\int_a^b\Gamma_{s,\varepsilon}(x)dx,\]
moreover the integral on the left hand side is taken in Lebesgue meaning. Since $F_s(x, y)$ is integrable, then the relations \eqref{GrindEQ__5_}, which are valid for any dissection, show:
\[\int_a^b\left(\int^{\ast}_{K\setminus \Sigma_{\varepsilon, s}}F_s(x, y)dy\right)dx=\int_a^b\Gamma_{s, \varepsilon}(x)dx=\iint_KF_s(x, y)dxdy.\]
From the previous relation and the last equality, we obtain, taking into account the equality \eqref{GrindEQ__4_}:
\[\iint_{K\setminus\Sigma_{\varepsilon}}f(x, y)dxdy=\int_a^b\Gamma_{\varepsilon}(x)dx,\]
moreover, at the right hand side the Riemann integral stands, that is the integral taken above in Lebesgue meaning exists in the Riemann sense, actually. To prove the last statement, fix some positive number $\varepsilon'$. It is possible find the natural number $s$ so that 
\[\sup|f(x, y)-f(x', y')|\le \lambda_s\]
for every pair of tuples $(x, y)$ and $(x', y')$, if they satisfy the inequality $|x-x'|+|y-y'|<\delta$. Letting $\lambda_s<\varepsilon'$ we dissect the rectangle $K$ into small rectangles every of which has a diameter $\le \delta$. Construct an integral sum:
\[\sum_{i=1}^n\left(\int^{\ast}F_s(\xi_i, y)dy\right)\Delta x_i.\]
From \eqref{GrindEQ__5_} is obvious that independently of the points $\xi_i$, this sum placed between lower and upper Darboux sums which differs each from other no more than $|\overline{D}_T(F_s)-\underline{D}_T(F_s)|\le\lambda_s|K|\le|K|\varepsilon'$, as it was shown above. Since $|K|$ is bounded this relation shows that the integral sum tends to the limit, that is the considered integral is a Riemann integral.

Define now some sequence of positive numbers  $(\varepsilon_n)$ tending to zero. The sequence of integrals $\Gamma_{x,\varepsilon_n}$ being written for every $\varepsilon_n$, are taken over increasing sequence of Jordan domains. So, from the boundedness of the function it follows the convergence of this sequence poinwisely, withrespect to $x$. From the relation  \eqref{GrindEQ__7_}, the relation below follows, by the theorem on bounded convergence:
\[\int _{a}^{b}\Gamma(x) dx=\lim_{n\to\infty}\int_a^b\Gamma_{x,\varepsilon_n}.\]
Now again, as above, it required to establish that the integral on the left hand side, taken over the segment $[a, b]$, exists in the Riemann sense and is equal to the initial double integral. 

Fix now, as ablve, some positive number $\varepsilon'$. There exists a natural number $s$ for which 
\[\sup|f(x, y)-f(x', y')|\le \lambda_s,\]
for every pair of pairs $(x, y)$ and $(x', y')$, if they are satisfying the condition $|x-x'|+|y-y'|<\delta$. Supposing that $\lambda_s<\varepsilon'$, consider some dissection of the seqment $[a, b]$ into small ones with maximal lengths $\le \delta$. Dissect also the segment $[c, d]$ by the same way. We suppose that the dissections are made by such manner that all quadrates of the covering $\Sigma_{\varepsilon_m,s}$, $\varepsilon_m<\varepsilon'$ also are dissected into small parts. Construct the integral sum:
\[\sum_{i=1}^p\left(\int^{\ast}f(\xi_i, y)dy\right)\Delta x_i.\]
Using the analog of the relations \eqref{GrindEQ__5_}, written for the function $f(x, y)$, we note that the introduced integral sum (without $\sup$) takes part in the chains of inequalities \eqref{GrindEQ__5_}. Therefore, exchanging upper integrals by integral sums below
\[\sum_{j=1}^m f(\xi_i,\eta_j)\Delta y_j,\]
we get the integral sum for double integral which is placed between upper and lower Darboux sums. The part of the integral sum where pairs of picked points $(\xi_i,\eta_j)$ fall into the quadrates of the covering $\Sigma_{\varepsilon_m,s}$, has a contribution $O(M\varepsilon')$. Remaining part of the integral sum coinsides with the integral sum of the function $F_s(x, y)$. For sufficiently small values $\varepsilon'$ this sum differs from the integral sum 
\[\sum_{i=1}^p\left(\int^{\ast}F_s(\xi_i, y)dy\right)\Delta x_i,\]
by a negligible error which is a quantity of the order $O((M+|K|)\varepsilon')$. From the relation \eqref{GrindEQ__6_} it follows that the limit does not depend on the points of dissection and the points which picked out arbitrarily. So, the last sum is near to the integral sum for the integral $\int _{a}^{b}\Gamma(x) dx$. Since $\varepsilon'$ is sufficiently small, these notes imply that the integral exists in the Riemann sense. Proof of Theorem 6 is finished.


\begin{thebibliography}{39}


\bibitem{1} R.~Courant, \textit{Cource of differential and integral calculus}, Moscow: Nauka, Vol. 2, 1970.(rus)

\bibitem{2} I. Sh.~Dzhabbarov. \textit{On multidimensional Tarry's problem for cubic polynomials}. Mathematical notes, \textbf{107} (5), 2020. pp. 657-673.

\bibitem{3} B. ~Gelbaum, J.~Olmsted, \textit{Counterexamples in analysis}. Moscow: Mir, 1967.(rus)

\bibitem{4} M.~K.~Grebencha, S.~I.~ Novoselov, \textit{Course of mathematical analysis}. Moscow: Visshaya shkola, Vol. 2, 1961.

\bibitem{5} I.~Sh.~Jabbarov, S.~A.~Meshaik, M.~M.~Ismailova, \textit{On the Number of Sheets of Coverings Defined by a System of Equations in n-Dimensional Spaces}. Chebishevskii sbornik, \textbf{24}(4), 2023.

\bibitem{6} S.~M.~Nikolskii, \textit{Course of mathematical analysis}, v.2. Moscow: Nauka, 1991.(rus)

\bibitem{7} G.~E.~Shilov, \textit{Mathematical analysis. Functions of several real variables}. Moscow: Nauka, 1972. (rus)

\bibitem{8} E.~A.~Zorich, \textit{Mathematical analysis}. Moscow: MCNMS, Vol. 2, 4th ed., 2002. (rus)

\bibitem{9}E.~Kamke, \textit{Handbook on the differential equations in partial derivatives of the first order.} Moscow: Nauka, 1966. (rus).

\bibitem{10} N. Bourbaki, \textit{Sketch on history of Mathematics.} Moscow: IIL, 1963. (rus). 

\end{thebibliography}
\end{document}